\documentclass[11pt]{article}
\usepackage{latexsym,amssymb,amsmath,graphics,cite,arydshln,amsfonts,multirow}
\usepackage{tikz}
\usepackage{fancyhdr}
\usepackage{lipsum}
\usepackage{lineno}
\usepackage{color}

\begin{document}
\newcommand{\qed}{\hphantom{.}\hfill $\Box$\medbreak}
\newcommand{\proof}{\noindent{\bf Proof \ }}
\newtheorem{Theorem}{Theorem}[section]
\newtheorem{Lemma}[Theorem]{Lemma}
\newtheorem{Corollary}[Theorem]{Corollary}
\newtheorem{Remark}[Theorem]{Remark}
\newtheorem{Example}[Theorem]{Example}
\newtheorem{Definition}[Theorem]{Definition}
\newtheorem{Construction}[Theorem]{Construction}


\begin{center}
{\Large\bf Optimal $2$-D $(n\times m,3,2,1)$-optical orthogonal codes and related equi-difference conflict avoiding codes  \footnote{Supported by NSFC under Grant $11471032$, and Fundamental Research Funds for the Central Universities under Grant $2016$JBM$071$, $2016$JBZ$012$ (T. Feng), NSFC under Grant 11401582, and NSFHB under Grant A2015507019 (L. Wang), NSFC under Grant 11771227, and Zhejiang Provincial Natural Science Foundation of China under Grant LY17A010008 (X. Wang).}
}

\vskip12pt

Tao Feng$^a$, Lidong Wang$^b$, and Xiaomiao Wang$^c$\\[2ex] {\footnotesize $^a$Department of Mathematics, Beijing Jiaotong University, Beijing 100044, P. R. China}
\\{\footnotesize $^b$Department of Basic Courses, Chinese People's Armed Police Force Academy, Langfang 065000, P. R. China}
\\{\footnotesize $^c$Department of Mathematics, Ningbo University, Ningbo 315211, P. R. China}\\
{\footnotesize
tfeng@bjtu.edu.cn, lidongwang@aliyun.com, wangxiaomiao@nbu.edu.cn}
\vskip12pt

\end{center}

\vskip12pt

\noindent {\bf Abstract:} This paper focuses on constructions for optimal $2$-D $(n\times m,3,2,1)$-optical orthogonal codes with $m\equiv 0\ ({\rm mod}\ 4)$. An upper bound on the size of such codes is established. It relies heavily on the size of optimal equi-difference $1$-D $(m,3,2,1)$-optical orthogonal codes, which is closely related to optimal equi-difference conflict avoiding codes with weight $3$. The exact number of codewords of an optimal $2$-D $(n\times m,3,2,1)$-optical orthogonal code is determined for $n=1,2$, $m\equiv 0 \pmod{4}$, and $n\equiv 0 \pmod{3}$, $m\equiv 8 \pmod{16}$ or $m\equiv 32 \pmod{64}$ or $m\equiv 4,20 \pmod{48}$.

\vskip12pt


\vskip12pt



\section{Introduction}

Let $n$, $m$, $k$, $\lambda_a$ and $\lambda_c$ be positive
integers. A {\em two-dimensional $(n\times m,k,\lambda_a,\lambda_c)$ optical
orthogonal code} (briefly $2$-D $(n\times m,k,\lambda_a,\lambda_c)$-OOC),
$\cal{C}$, is a family of $n\times m$ $(0, 1)$-matrices (called {\em
codewords}) of Hamming weight $k$ satisfying the following two properties:
\begin{center}
\begin{enumerate}
\item[(1)] the autocorrelation property: for each matrix
${\mathbf{A}}=(a_{ij})_{n\times m}\in\cal{C}$ and each
integer $r$, $r\not\equiv 0\ ({\rm mod}\ m)$,
$$\sum_{i=0}^{n-1}\sum_{j=0}^{m-1}a_{ij}a_{i,j+r}\leq\lambda_a;$$

\item[(2)] the cross-correlation property: for each matrix
${\mathbf{A}}=(a_{ij})_{n\times m}\in\cal{C}$, ${\mathbf{B}}=(b_{ij})_{n\times m}\in\cal{C}$ with ${\mathbf{A}}\neq {\mathbf{B}}$, and each integer $r$,
$$\sum_{i=0}^{n-1}\sum_{j=0}^{m-1}a_{ij}b_{i,j+r}\leq\lambda_c.$$
\end{enumerate}
\end{center}
where the arithmetic $j+r$ is reduced modulo $m$. When $n=1$, a two-dimensional $(1\times m,k,\lambda_a,\lambda_c)$ optical orthogonal code is said to be a \emph{one-dimensional $(m,k,\lambda_a,\lambda_c)$-optical orthogonal code}, denoted by a $1$-D $(m,k,\lambda_a,\lambda_c)$-OOC.

Optical orthogonal codes are widely used as spreading codes in optical fiber networks. $1$-D OOC was first investigated systematically by Chung Salehi and Wei in \cite{csw}. $1$-D OOCs have a drawback which requires a large chip rate. To overcome it, $2$-D OOCs were proposed in \cite{yk}, which spreads in both time and wavelength so that the chip rate requirement can be substantially reduced.

The number of codewords of a $2$-D OOC is called its {\em size}. For fixed $n$, $m$, $k$, $\lambda_a$ and $\lambda_c$, the largest size among all $2$-D
$(n\times m,k,\lambda_a,\lambda_c)$-OOCs is denoted by $\Phi(n\times m,k,\lambda_a,\lambda_c)$. A $2$-D $(n\times m,k,\lambda_a,\lambda_c)$-OOC with $\Phi(n\times m,k,\lambda_a,\lambda_c)$ codewords is said to be {\em
optimal}. Naturally, a $1$-D $(m,k,\lambda_a,\lambda_c)$-OOC is said to be {\em optimal} if it contains $\Phi(1\times m,k,\lambda_a,\lambda_c)$ codewords.

When $\lambda_a=\lambda_c=\lambda$, various 2-D OOCs or 1-D OOCs were constructed based on algebraic and combinatorial methods (see \cite{ab1,am1,am2,be,b2,b18,bp2,cc1,cfm,cm,cws,fc,fcj,fm,gy,hc,hc1,mms,ogkeb,wc,wsy,wy,y} and the references therein). Instead, very little has been done on optimal OOCs with $\lambda_a\neq\lambda_c$. Yang and Fuja \cite{yf} showed that the auto- and cross-correlation properties are used for synchronization and user identification, respectively, and in some circumstances only with good cross-correlation one can deal with both synchronization and user identification. This motivates the study of OOCs with better cross-correlation than auto-correlation. See \cite{bmp,bpw,mb,wc1} for example, in which the cases of $n=1$, $k\in\{4,5\}$, $\lambda_a=2$ and $\lambda_c=1$ are considered.

When $\lambda_a=k$ and $\lambda_c=1$, a $1$-D $(m,k,k,1)$-OOC is also called a {\em conflict-avoiding code}, denoted by a CAC$(m,k)$, which can be viewed as a $1$-D $(m,k,1)$-OOC without the constraint of the auto-correlation property. A CAC finds its application on a multiple-access collision channel without feedback (cf. \cite{gv,l}).

When $m$ is even, optimal CAC$(m,3)$s have been discussed thoroughly in \cite{lt,jimbo,mfu,fu}. We summarize the results for later use.

\begin{Theorem}{\rm \cite{lt,jimbo,mfu,fu}}
\label{cac} Let $m\equiv 0\ ({\rm mod}\ 2)$.
The size of an optimal CAC$(m,3)$ $($i.e., an optimal $1$-D $(m,3,3,1)$-OOC$)$ is
$$\Phi(1\times m,3,3,1)=\left\{
\begin{array}{lll}
(m-2)/4, & m\equiv 2\ ({\rm mod}\ 4),\\
\lfloor (7m+16)/32\rfloor, & m\equiv 0\ ({\rm mod}\ 24) \ {\rm and}\ m\neq 48,\\
\lfloor (7m+4)/32\rfloor, & m\equiv 4,20\ ({\rm mod}\ 24),\\
\lfloor 7m/32\rfloor, & m\equiv 8,16\ ({\rm mod}\ 24)\ {\rm and}\ m\neq 64,\\
\lfloor (7m+20)/32\rfloor, & m\equiv 12\ ({\rm mod}\ 24),\\
\end{array}
\right.
$$
with the exception of $\Phi(1\times 48,3,3,1)=10$ and
$\Phi(1\times 64,3,3,1)=13.$
\end{Theorem}

However, there are few results on optimal $2$-D $(n\times m,k,2,1)$-OOCs when $n\neq 1$ in the literature. The only known results for $k=3$ is from \cite{fwwz,wcf}, which determined the size of an optimal $2$-D $(n\times m,3,2,1)$-OOCs with $m\equiv 2\ ({\rm mod}\ 4)$.
This paper continues the work in \cite{fwwz}, and we are concerned about optimal $2$-D $(n\times m,3,2,1)$-OOCs with $m\equiv 0\ ({\rm mod}\ 4)$.

In Section $2$, an equivalent description of $2$-D $(n\times m,k,\lambda_a,1)$-OOCs is given by using set-theoretic notation. Section $3$ is devoted to presenting an upper bound on the size of an optimal $2$-D $(n\times m,3,2,1)$-OOC with $m\equiv 0\ ({\rm mod}\ 4)$. We will see that the upper bound relies heavily on the size of an optimal equi-difference $1$-D $(m,3,2,1)$-OOC (see Lemma \ref{0mod4}). So we focus our attention on constructions for such kind of 1-D OOCs in Section 4. Interestingly, equi-difference $1$-D $(m,3,2,1)$-OOCs are closely related to equi-difference CAC$(m,3)$s (see Section 4.2). The latter have been investigated independently in recent times \cite{fls,lmsj,mzs,m0,wf}.

In Section 5, direct and recursive constructions for optimal $2$-D $(n\times m,3,2,1)$-OOCs are given. We shall point out why it seems to be difficult to find effective recursive constructions for optimal $2$-D $(n\times m,3,2,1)$-OOCs (see Remark \ref{remark:difficult}).

Throughout this paper, assume that $I_n=\{0,1,\ldots,n-1\}$ and denote by $Z_m$ the additive group of integers modulo $m$. For a unit $a\in Z_m\setminus\{0\}$, the \emph{multiplicative order} of $a$, denoted by ord$_m(a)$, is the smallest positive integer $l$ such that $a^l\equiv 1\ ({\rm mod}\ m)$. As the main result of this paper, we are to prove the following theorem.

\begin{Theorem}\label{main result} The size of an optimal $2$-D $(n\times m,3,2,1)$-OOC is
$$\Phi(n\times m,3,2,1)=\left\{
\begin{array}{lll}
\lfloor 7m/32\rfloor, &  n=1, m\equiv 0\ ({\rm mod}\ 8) {\rm\ and\ } m\neq64;\\
\lfloor (7m+4)/32\rfloor, &  n=1 {\rm\ and\ } m\equiv 4\ ({\rm mod}\ 8);\\
13, &  (n,m)=(1,64);\\
3m/4, &  n=2, m\equiv 0\ ({\rm mod}\ 4) {\rm\ and\ } m>4;\\
2, &  (n,m)=(2,4);\\
n(8nm+3m-8)/48, &  n\equiv 0\ ({\rm mod}\ 3), n\neq 6,9, {\rm\ and\ } \\
& m\equiv 8\ ({\rm mod}\ 16);\\
n(32nm+11m-32)/192, &  n\equiv 0\ ({\rm mod}\ 3), n\neq 6,9, {\rm\ and\ } \\
& m\equiv 32\ ({\rm mod}\ 64);\\
6, &  (n,m)=(3,4);\\
n(8nm+3m+4)/48, &  n\equiv 0\ ({\rm mod}\ 3),  n\neq 6,9, m>4, \\
& m\equiv 4,20\ ({\rm mod}\ 48), {\rm\ and\ }  m/4\in S\\
\end{array}
\right.
$$
codewords, where $S$ is the set of positive integers such that for any $s\in S$, it holds that $s\equiv 1,5\ ({\rm mod}\ 12)$, and every prime divisor $p$ of $s$ satisfies $p\equiv 5\ ({\rm mod}\ 8)$, or $p\equiv 1\ ({\rm mod}\ 8)$ and $4|ord_p(2)$.
\end{Theorem}

\section{Preliminaries}

A convenient way of viewing optical orthogonal codes is from a set-theoretic perspective.

Let $\cal C$ be a $2$-D $(n\times m,k,\lambda_a,\lambda_c)$-OOC. For each $n\times m$
$(0,1)$-matrix $M\in\cal C$, whose rows are indexed by $I_n$ and
columns are indexed by $Z_m$. Construct a $k$-subset $B_M$ of $I_n\times Z_m$ such that $(i,j)\in B_M$ if and only if $M$'s $(i,j)$ cell equals $1$. Then $\{B_M: M\in {\cal C}\}$ is a set-theoretic representation of the $2$-D $(n\times m,k,\lambda_a,\lambda_c)$-OOC. Conversely, let $\cal B$ be a set of $k$-subsets of $I_n\times Z_m$. ${\cal B}$ constitutes a $2$-D $(n\times m,k,\lambda_a,\lambda_c)$-OOC if the following two conditions are satisfied:
\begin{enumerate}
\item[($1'$)] the autocorrelation property: $|B\cap(B+s)|\leq\lambda_a$ for any $B\in{\cal B}$ and any integer $s$, $s\not\equiv 0\ ({\rm mod}\ m)$;
\item[($2'$)] the cross-correlation property: $|A\cap(B+s)|\leq\lambda_c$ for any $A,B\in{\cal B}$ with $A\neq B$ and any integer $s$,
\end{enumerate}
where $B+s=\{(i,x+s)\ \pmod{(-,m)}:(i,x)\in B\}$.

It is not convenient to check the autocorrelation and cross-correlation property of a set $\cal B$ of $k$-subsets of $I_n\times Z_m$ via Conditions $(1')$ and $(2')$. However, when $\lambda_c=1$, one can use the pure and mixed difference method to describe a $2$-D $(n\times m,k,\lambda_a,1)$-OOC.

For $(i,x),(i,y)\in I_n\times Z_m$ with $x\neq y$, the difference $x-y$ (mod $m$) is called a {\em pure $(i,i)$-difference}. For $(i,x),(j,y)\in I_n\times Z_m$ with $i\neq j$, the difference $x-y$ (mod $m$) is called a {\em mixed $(i,j)$-difference}. Let $B$ be a $k$-subset of $I_n\times Z_m$. Given $i,j\in I_n$, define a multi-set
$$\Delta_{ij}(B)=\{x-y\ (mod\ m): (i,x),(j,y)\in B, (i,x)\neq(j,y)\}.$$ When $i=j$, $\Delta_{ii}(B)$ is the multi-set of all pure $(i,i)$-differences of $B$. When $i\neq j$, $\Delta_{ij}(B)$ is the multi-set of all mixed $(i,j)$-differences of $B$. Note that $\Delta_{ij}(B)$ is empty if $i$ or $j$ does not occur as the first component of the elements of $B$.

Let $\cal B$ be a set of $k$-subsets of $I_n\times Z_m$ and $B\in {\cal B}$. Let $\lambda(B)$ denote the maximum multiplicity of elements in the multi-set $\bigcup_{i\in I_n}\Delta_{ii}(B)$. Then ${\cal B}$ constitutes a $2$-D $(n\times m,k,\lambda_a,1)$-OOC if the following two conditions are satisfied:

\begin{enumerate}
\item[$(1'')$] the autocorrelation property: $\lambda(B)\leq\lambda_a$ for any $B\in{\cal B}$;
\item[$(2'')$] the cross-correlation property: $\Delta_{ij}(A)\cap \Delta_{ij}(B)=\emptyset$ for any $A,B\in{\cal B}$ with $A\neq B$ and any $i,j\in I_n$ ($i$ may be equal to $j$).
\end{enumerate}

\noindent The interested reader is referred to \cite{hc1} for details on the equivalence of $(1')$ and $(1'')$.

In the remainder of this paper, we always use the set-theoretic language to describe $2$-D OOCs.

\section{Upper bound on the size of $2$-D $(n\times m,3,2,1)$-OOCs}

In a $2$-D $(n\times m,3,2,1)$-OOC, each codeword is of the form $\{(i_1,x),(i_2,y),(i_3,z)\}$, where $i_1,i_2,i_3\in I_n$ and $x,y,z\in Z_m$. All codewords can be divided into the following three types:

\begin{itemize}
\item Type $1$: $i_1=i_2=i_3$;
\item Type $2$: $i_1=i_2\neq i_3$;
\item Type $3$: $i_1$, $i_2$, $i_3$ are pairwise distinct.
\end{itemize}

\noindent Let $\alpha$, $\beta$, $\gamma$ denote the numbers of
codewords of Type $1$, $2$, $3$ in a $2$-D $(n\times m,3,2,1)$-OOC, respectively.

For Type $1$, the codewords can be classified further according
to the second coordinates. Take any codeword $\{(i_1,x),(i_1,y),(i_1,z)\}$ of Type $1$ and consider its derived set $X=\{x,y,z\}$ of the second coordinates. Define the list of
differences of $X$ by $\Delta X=\{b-a\ ({\rm mod}\ m):\ a,b\in X,a\neq b\}$ as a
multi-set, and define the \emph{support} of $\Delta X$, denoted by
${\rm supp}(\Delta X)$, as the set of underlying elements in $\Delta X$. Define the {\em orbit} of $X$ under $Z_m$ by ${\rm Orb}(X)=\{\{x+i\pmod{m}:x\in X\}:i\in Z_m\}$.
By Lemma $2.2$ in \cite{mb}, we have
$$|{\rm supp}(\Delta X)|=\left\{
\begin{array}{lll}
2, & X\in {\rm Orb}(\{0,m/3,2m/3\}),\\
3, & X\in {\rm Orb}(\{0,m/4,m/2\}),\\
4, & X\in {\rm Orb}(\{0,a,2a\}) {\rm \ except\ for\ the\ cases\ of\ } |{\rm supp}(\Delta
X)|=2,3,\\
5, & X\in {\rm Orb}(\{0,a,m/2\}) {\rm \ except\ for\ the\ case\ of\ }
|{\rm supp}(\Delta X)|=3,\\
6, & X\in {\rm Orb}(\{0,a,b\}) {\rm \ except\ for\ the\ cases\ of\ }
|{\rm supp}(\Delta X)|=2,3,4,5.\\
\end{array}
\right.
$$
If $|{\rm supp}(\Delta X)|=2$, then $X\in {\rm Orb}(\{0,m/3,2m/3\})$, which implies that $m/3$ occurs three times as a pure $(i_1,i_1)$-difference. It contradicts with the autocorrelation parameter $\lambda_a=2$. Thus $|{\rm supp}(\Delta X)|=3,4,5$ or $6$. Let $\alpha_3$, $\alpha_4$, $\alpha_5$, $\alpha_6$ denote the numbers of codewords of Type $1$ in a $2$-D $(n\times m,3,2,1)$-OOC such that each derived set $X$ of these codewords satisfies $|{\rm supp}(\Delta X)|=3,4,5,6$, respectively. Then $\alpha_3+\alpha_4+\alpha_5+\alpha_6=\alpha$.

For Type $2$, take any codeword $\{(i_1,x),(i_1,y),(i_2,z)\}$ with $i_1\neq i_2$ and consider its partial derived set $Y=\{x,y\}$ of the second coordinates. Let $\beta_1$ denote the number of codewords of Type $2$ in a $2$-D $(n\times m,3,2,1)$-OOC such that each partial derived set $Y$ of these codewords satisfies $y-x\equiv m/2\ ({\rm mod}\ m)$. Denote by $\beta_2$ the number of the remaining codewords of Type $2$ in the $2$-D OOC. Then $\beta_1+\beta_2=\beta$.

\subsection{General upper bound}

We need a new concept. A $1$-D $(m,3,2,1)$-OOC is said to be \emph{equi-difference} if each of its codewords is of the form $X=\{0,a,2a\}$, i.e., $|supp(\Delta X)|=3$ or $4$. Let $\Psi^e(m,3,2,1)$ denote the largest size of codes among all equi-difference $1$-D $(m,3,2,1)$-OOCs for given $m$. An equi-difference $1$-D $(m,3,2,1)$-OOC is said to be \emph{optimal} if it contains $\Psi^e(m,3,2,1)$ codewords.

\begin{Lemma} \label{alpha} $(1)$ $\alpha_3+\alpha_4\leq n\Psi^e(m,3,2,1)$.

$(2)$ $\alpha_3+\alpha_5+\beta_1\leq n$.
\end{Lemma}

\proof $(1)$ Examine codewords of the form $\{(i,0),(i,a),(i,2a)\}$ with $X=\{0,a,2a\}$ and $|{\rm supp}(\Delta X)|=3,4$. The number of such kind of codewords is not more than $\Psi^e(m,3,2,1)$ for each $i\in I_n$. So $\alpha_3+\alpha_4\leq n\Psi^e(m,3,2,1)$.

(2) For each $i\in I_n$, there is at most one codeword that admits $m/2$ as a pure $(i,i)$-difference. So $\alpha_3+\alpha_5+\beta_1\leq n$. \qed

\begin{Lemma} \label{0mod4}
$\Phi(n\times m,3,2,1)\leq\left\{
\begin{array}{lll}
\lfloor n(nm+2\Psi^e(m,3,2,1))/6\rfloor, & {\rm if\ } m\equiv 0\ ({\rm mod}\ 2),\\
\lfloor n(nm+2\Psi^e(m,3,2,1)-1)/6\rfloor, & {\rm if\ } m\equiv 1\ ({\rm mod}\ 2).\\
\end{array}
\right.$
\end{Lemma}

\proof In a $2$-D $(n\times m,3,2,1)$-OOC, given $i\in I_n$, there are at most $m-1$ different pure $(i,i)$-differences; and given $i,j\in I_n$ with $i\neq j$, there are at most $m$ different mixed $(i,j)$-differences. Thus the total numbers of different pure differences and mixed differences in a $2$-D $(n\times m,3,2,1)$-OOC are at most $n(m-1)$ and $n(n-1)m$, respectively. Pure differences are from Type $1$ and a part of Type $2$, while mixed differences are from Type $3$ and the other part of Type $2$. So we have
\begin{align}
\label{eqn1} 3\alpha_3+4\alpha_4+5\alpha_5+6\alpha_6+\beta_1+2\beta_2 & \leq n(m-1),  \\
\label{eqn2} 4\beta+6\gamma & \leq n(n-1)m.
\end{align}
By Lemma \ref{alpha}(2),
\begin{align}
\label{eqn3} \alpha_3+\alpha_5+\beta_1 & \leq n.
\end{align}
Note that $\alpha_3+\alpha_4+\alpha_5+\alpha_6=\alpha$ and $\beta_1+\beta_2=\beta$.
By (\ref{eqn1})+(\ref{eqn2})+(\ref{eqn3}), we have $6(\alpha+\beta+\gamma)-2(\alpha_3+\alpha_4)\leq
n^2m$. By Lemma \ref{alpha}(1), $\alpha_3+\alpha_4\leq n\Psi^e(m,3,2,1)$. It follows that
$\alpha+\beta+\gamma\leq
(n^2m+2n\Psi^e(m,3,2,1))/6.$
Therefore, $\Phi(n\times m,3,2,1)\leq\lfloor n(nm+2\Psi^e(m,3,2,1))/6\rfloor$.

Furthermore, when $m$ is odd, $\alpha_3=\alpha_5=\beta_1=0$. So $\alpha_4+\alpha_6=\alpha$ and $\beta_2=\beta$. Then by (\ref{eqn1})+(\ref{eqn2}), we have $6(\alpha+\beta+\gamma)-2\alpha_4\leq n^2m-n$. By Lemma \ref{alpha}(1), $\alpha_4\leq n\Psi^e(m,3,2,1)$. Thus $\Phi(n\times m,3,2,1)\leq\lfloor n(nm+2\Psi^e(m,3,2,1)-1)/6\rfloor$ for any odd integer $m$. \qed

\begin{Remark} \label{remark}
Examining the proof of Lemma $\ref{0mod4}$, we have that if
$$\Phi(n\times m,3,2,1)=\left\{
\begin{array}{lll}
n(nm+2\Psi^e(m,3,2,1))/6 & {\rm if\ } m\equiv 0\ ({\rm mod}\ 2),\\
n(nm+2\Psi^e(m,3,2,1)-1)/6, & {\rm if\ } m\equiv 1\ ({\rm mod}\ 2),
\end{array}
\right.$$
then $\alpha_3+\alpha_4= n\Psi^e(m,3,2,1)$ from Lemma $\ref{alpha}(1)$, which means that in such cases, any optimal $2$-D $(n\times m,3,2,1)$-OOC must contain $n$ optimal equi-difference $1$-D $(m,3,2,1)$-OOCs as subcodes.
\end{Remark}

\begin{Lemma}\label{3*4-bound}
$\Phi(3\times 4,3,2,1)\leq 6$.
\end{Lemma}

\proof An optimal equi-difference $1$-D $(4,3,2,1)$-OOC defined on $Z_4$ contains only one codeword $\{0,1,2\}$, so $\Psi^e(4,3,2,1)=1$. Then by Lemma \ref{0mod4}, $\Phi(3\times 4,3,2,1)\leq 7$. Assume that $\Phi(3\times 4,3,2,1)=7$. By Remark \ref{remark}, $\alpha_3+\alpha_4=3$. Since $\alpha_4=\alpha_5=\alpha_6=0$ for any $2$-D $(3\times 4,3,2,1)$-OOC, we have $\alpha_3=3$. By Formula (\ref{eqn1}), $3\alpha_3+4\alpha_4+5\alpha_5+6\alpha_6+\beta_1+2\beta_2 \leq 9$, so $\beta=0$, which yields $\gamma=4$. Write the $4$ codewords of Type 3 as $\{(0,0),(1,a_i),(2,b_i)\}$, $i=1,2,3,4$.
Clearly, $\bigcup_{i=1}^4 \{a_i\}=\bigcup_{i=1}^4 \{b_i\}=\bigcup_{i=1}^4 \{b_i-a_i\ ({\rm mod}\ 4)\}=Z_4$. Thus $\sum_{i=1}^4 (b_i-a_i)=\sum_{i=1}^4 b_i-\sum_{i=1}^4 a_i=0$ and $\sum_{i=1}^4 (b_i-a_i)\equiv 0+1+2+3\ ({\rm mod}\ 4)$, a contradiction. \qed

\subsection{Improved upper bound for $n=2$}

\begin{Lemma}\label{n=2bound}
$\Phi(2\times m,3,2,1)\leq\left\{
\begin{array}{lll}
\lfloor 3m/4 \rfloor, & {\rm if\ } m\equiv 0\ ({\rm mod}\ 2),\\
\lfloor (3m-2)/4\rfloor, & {\rm if\ } m\equiv 1\ ({\rm mod}\ 2).\\
\end{array}
\right.$
\end{Lemma}

\proof Formulas (\ref{eqn1})-(\ref{eqn3}) in Lemma \ref{0mod4} still hold when $n=2$. Note that $\gamma=0$ when $n=2$. We rewrite these formulas as follows
\begin{align}
\label{eqn4} 3\alpha_3+4\alpha_4+5\alpha_5+6\alpha_6+\beta_1+2\beta_2 & \leq 2(m-1), \\
\label{eqn5} 2\beta & \leq m, \\
\label{eqn6} \alpha_3+\alpha_5+\beta_1 & \leq 2.
\end{align}
By (\ref{eqn4})+(\ref{eqn5})+(\ref{eqn6}), we have
$4(\alpha+\beta)+2(\alpha_5+\alpha_6)\leq3m$. Due to
$\alpha_5,\alpha_6\geq 0$, we have $\alpha+\beta\leq 3m/4$. Hence,
$\Phi(2\times m,3,2,1)\leq \lfloor3m/4\rfloor$.

Furthermore, when $m$ is odd, $\alpha_3=\alpha_5=\beta_1=0$. Then by (\ref{eqn4})+(\ref{eqn5}), we have $4(\alpha+\beta)+2\alpha_6\leq 3m-2$. Thus $\Phi(n\times m,3,2,1)\leq\lfloor (3m-2)/4\rfloor$ for any odd integer $m$. \qed

\begin{Lemma}\label{2*4}
$\Phi(2\times 4,3,2,1)\leq 2$.
\end{Lemma}

\proof By Lemma \ref{n=2bound}, $\Phi(2\times 4,3,2,1)\leq 3$. Assume that $\Phi(2\times 4,3,2,1)=3$. Since $\alpha_4=\alpha_5=\alpha_6=0$ for any $2$-D $(2\times 4,3,2,1)$-OOC, we rewrite Formulas (\ref{eqn4})-(\ref{eqn6}) as follows
\begin{align}
\label{eqn7} 3\alpha_3+\beta_1+2\beta_2 & \leq 6,\\
\label{eqn8} \beta & \leq 2,  \\
\label{eqn9} \alpha_3+\beta_1 & \leq 2.
\end{align}
$\Phi(2\times 4,3,2,1)=3$ yields $\alpha+\beta=3$, so $\alpha=\alpha_3\geq 1$ by (\ref{eqn8}). If $\alpha_3=2$, then $\beta=0$ by (\ref{eqn7}), which implies $\alpha+\beta=2$, a contradiction. So $\alpha_3=1$ and $\beta=2$. Then $\beta_1\leq 1$ by (\ref{eqn9}) and $\beta_2\geq 1$ by $\beta=2$. It follows that $\beta_1=\beta_2=1$ by (\ref{eqn7}). W.l.o.g., let the codeword such that $\alpha_3=1$ be $\{(0,0),(0,1),(0,2)\}$, and the codewords such that $\beta_1=1$ and $\beta_2=1$ are $\{(1,0),(1,2),(0,x)\}$ and $\{(1,0),(1,1),(0,y)\}$, respectively, for some $x,y\in Z_4$. Examining the mixed $(0,1)$-differences, we obtain $\{x,x-2,y,y-1\}\equiv\{0,1,2,3\} \pmod{4}$. It is readily checked that such $x$ and $y$ do not exist, a contradiction. \qed

\subsection{Improved upper bound for $n=1$ and $m\equiv 0 \pmod{4}$}

To present an improved upper bound for $\Phi(1\times m,3,2,1)$, we here review the linear programming approach formulated by Jimbo et al. \cite{jimbo}.

For any codeword $X$ in a $1$-D $(m,3,2,1)$-OOC with $m\equiv 0 \pmod{4}$, since the elements of $\Delta(X)$ are symmetric with respect to $m/2$, it suffices to consider the halved difference set $$\Delta_2(X)=\{i:\ i\in \Delta(X),1\leq i\leq m/2\}$$ instead of
$\Delta(X)$. Note that $\Delta(X)$ is a multi-set, but $\Delta_2(X)$ is not.

Now partition the positive integers not exceeding $m/2$ into the following three subsets:
\begin{center}
$O=\{i:\ i\equiv 1 \pmod{2}, 1\leq i\leq m/2\}$,

$E=\{i:\ i\equiv 2 \pmod{4}, 1\leq i\leq m/2\}$,

$D=\{i:\ i\equiv 0 \pmod{4}, 1\leq i\leq m/2\}$.
\end{center}
The integers in $O$ are odd, those in $E$ are said to be \emph{singly even} and those in $D$ are said to be \emph{doubly even}. It follows that any codeword of a $1$-D $(m,3,2,1)$-OOC can be categorized into the following two lemmas according to the halved difference set produced from it.

\begin{Lemma}\label{center} {\rm \cite{jimbo}}
Let $m\equiv 0 \pmod{4}$. Any codeword $X$ of the form $\{0,i,2i\}$
satisfying $\Delta_2(X)=\{i,j\}$, where $j=2i$ if $1\leq i\leq m/4$, and
$j=m-2i$ if $m/4<i<m/2$ and $i\neq m/3$, belongs to one of the
following three types:
\begin{enumerate}
\item[$(i)$] $i\in O$ and $j\in E$,
\item[$(ii)$] $i\in E$ and $j\in D$,
\item[$(iii)$] $i,j\in D$.
\end{enumerate}
\end{Lemma}

\begin{Lemma}\label{no-center} {\rm \cite{jimbo}}
Let $m\equiv 0 \pmod{4}$. Any codeword $X$ satisfying $\Delta_2(X)=\{i,j,k\}$ belongs to one of the following four types:
\begin{enumerate}
\item[$(iv)$] two of $i,j$ and $k$ are in $O$ and one is in $E$,
\item[$(v)$] two of $i,j$ and $k$ are in $O$ and one is in $D$,
\item[$(vi)$] two of $i,j$ and $k$ are in $E$ and one is in $D$,
\item[$(vii)$] $i,j,k\in D$.
\end{enumerate}
\end{Lemma}

Take a $1$-D $(m,3,2,1)$-OOC $\cal C$. Let $C_o$, $C_e$ and $C_d$ denote the sets of codewords in $\cal C$ of Types $(i)$, $(ii)$ and $(iii)$, respectively, and $N_{oe}$, $N_{od}$, $N_e$ and $N_d$ denote the sets of codewords in $\cal C$ of Types $(iv)$, $(v)$, $(vi)$ and $(vii)$, respectively. Note that any codeword $X\in\mathcal{C}$ with $|{\rm supp}(\Delta X)|=3$ or $4$ satisfies Lemma \ref{center}, while any codeword $X$ with $|{\rm supp}(\Delta X)|=5$ or $6$ satisfies Lemma \ref{no-center}. Then $$|\mathcal{C}|=|C_o|+|C_e|+|C_d|+|N_{oe}|+|N_{od}|+|N_e|+|N_d|.$$

\begin{Lemma}\label{bound1}
$$\Phi(1\times m,3,2,1)\leq\left\{
\begin{array}{lll}
\lfloor 7m/32\rfloor, & {\rm if\ } m\equiv\ 0\ ({\rm mod}\ 8),\\
\lfloor (7m+4)/32\rfloor, & {\rm if\ } m\equiv\ 4\ ({\rm mod}\ 8).\\
\end{array}
\right.
$$
\end{Lemma}

\proof A $1$-D $(m,3,2,1)$-OOC with $m\equiv 0 \pmod{4}$ contributes at most $m/4$ different odd differences that are not more than $m/2$, $\lceil m/8\rceil$ different singly even differences that are not more than $m/2$, and $\lfloor m/8\rfloor$ different doubly even differences that are not more than $m/2$. It follows that
\begin{align}
\label{eqn10}
|C_o|+2|N_{oe}|+2|N_{od}|\leq m/4,\\
\label{eqn11}
|C_o|+|C_e|+|N_{oe}|+2|N_e|\leq\lceil m/8\rceil,\\
\label{eqn12}
|C_e|+2|C_d|+|N_{od}|+|N_e|+3|N_d|\leq\lfloor m/8\rfloor.
\end{align}
By (\ref{eqn10})+3(\ref{eqn11})+2(\ref{eqn12}), we have
$$4|\mathcal{C}|+|C_e|+|N_{oe}|+4|N_e|+2|N_d|\leq\left\{
\begin{array}{lll}
7m/8, & {\rm if\ } m\equiv\ 0\ ({\rm mod}\ 8), \\
(7m+4)/8, & {\rm if\ } m\equiv\ 4\ ({\rm mod}\ 8), \\
\end{array}
\right.
$$
where $|\mathcal{C}|=|C_o|+|C_e|+|C_d|+|N_{oe}|+|N_{od}|+|N_e|+|N_d|$ is the total number of codewords. Hence, $$|\mathcal{C}|\leq\left\{
\begin{array}{lll}
\lfloor 7m/32\rfloor, & {\rm if\ } m\equiv\ 0\ ({\rm mod}\ 8),\\
\lfloor (7m+4)/32\rfloor, & {\rm if\ } m\equiv\ 4\ ({\rm mod}\ 8).\\
\end{array}
\right.
$$ This completes the proof. \qed

\section{Equi-difference $1$-D $(m,3,2,1)$-OOCs}

By Lemma \ref{0mod4}, it is important to determine the exact value of $\Psi^e(m,3,2,1)$. Clearly, $\Psi^e(m,3,2,1)\leq\Phi(1\times m,3,2,1)$. A better upper bound can be shown in the following lemma.

\begin{Lemma} \label{bounf equi-diff}
$$\Psi^e(m,3,2,1)\leq\left\{
\begin{array}{ll}
\lfloor (m-1)/4\rfloor, &  {\rm if\ } m\not\equiv 0\ ({\rm mod}\ 4);\\
\lceil m/8\rceil+\Psi^e(m/4,3,2,1), &  {\rm if\ } m\equiv 0\ ({\rm mod}\ 4).\\
\end{array}
\right.
$$
\end{Lemma}

\proof Let $\cal C$ be an equi-difference $1$-D $(m,3,2,1)$-OOC. When $m\not\equiv 0\pmod{4}$, for any codeword $X\in{\cal C}$, $|{\rm supp}(\Delta X)|=4$, so $\Psi^e(m,3,2,1)\leq \lfloor (m-1)/4\rfloor$.

When $m\equiv 0\pmod{4}$, recall that $C_o$, $C_e$ and $C_d$ denote the sets of codewords in $\cal C$ of Types $(i)$, $(ii)$ and $(iii)$ (see Lemma \ref{center}), respectively. Then $|{\cal C}|=|C_o|+|C_e|+|C_d|$. A $1$-D $(m,3,2,1)$-OOC with $m\equiv 0 \pmod{4}$ contributes at most $\lceil m/8\rceil$ different singly even differences that are not exceeding $m/2$, which gives $|C_o|+|C_e|\leq \lceil m/8\rceil$. Observing that $|C_d|\leq \Psi^e(m/4,3,2,1)$, we obtain $|{\cal C}|\leq \lceil m/8\rceil+\Psi^e(m/4,3,2,1)$. \qed

For an equi-difference $1$-D $(m,3,2,1)$-OOC, $\cal B$, on $Z_m$, define $\Delta({\cal B})=\bigcup_{B\in{\cal B}}\Delta(B)$ to be a multi-set of differences, and the \emph{support} of $\Delta({\cal B})$, written as ${\rm supp}(\Delta({\cal B}))$, to be the set of underlying elements in $\Delta({\cal B})$. The {\em difference leave} of $\cal B$ is a set that consists of all nonzero elements of $Z_m$ not covered by ${\rm supp}(\Delta({\cal B}))$. Let $m\equiv 0\ ({\rm mod}\ g)$ and $H$ be the subgroup of order $g$ in $Z_m$, i.e., $H=\{0,m/g,\ldots,(g-1)m/g\}$. If ${\rm supp}(\Delta({\cal B}))\subseteq Z_m\setminus H$, then $\cal B$ is said to be a {\em $g$-regular} equi-difference $1$-D $(m,3,2,1)$-OOC.

\begin{Lemma}\label{lem:equi-2 mod 4}
There is an optimal equi-difference $1$-D $(m,3,2,1)$-OOC with $\Psi^e(m,3,2,1)$ $=(m-2)/4$ codewords for any $m\equiv 2\ ({\rm mod}\ 4)$, whose difference leave is $\{m/2\}$.
\end{Lemma}

\proof It is readily checked that $\{\{0,i,2i\}$, $i=1,3,\ldots,m/2-2\}$ forms a $2$-regular equi-difference $1$-D $(m,3,2,1)$-OOC with $(m-2)/4$ codewords, whose difference leave is $\{m/2\}$. By Lemma $\ref{bounf equi-diff}$, it is optimal. \qed

\subsection{A recursive construction}

Let $A$ be a set of integers and $w$ be an integer. Write $w\cdot A=\{wa:a\in A\}$. The following construction is straightforward by the definition of $g$-regular equi-difference $1$-D OOCs.

\begin{Construction} \label{filling g-regular} Suppose that there exist
\begin{enumerate}
\item[$(1)$] a $g$-regular equi-difference $1$-D $(m,3,2,1)$-OOC with $b_1$ codewords, whose difference leave is $L_1$ $($defined on $Z_m)$;
\item[$(2)$] an equi-difference $1$-D $(g,3,2,1)$-OOC with $b_2$ codewords, whose difference leave is $L_2$ $($defined on $Z_g)$.
\end{enumerate}
Then there exists an equi-difference $1$-D $(m,3,2,1)$-OOC with $b_1+b_2$ codewords, whose difference leave is $L_1\cup ((m/g)\cdot L_2)$ $($defined on $Z_m)$.
\end{Construction}

Assume that $[a,b]$ denotes the set of integers $n$ such that $a\leq n\leq b$, and $[a,b]_o$ denotes the set of odd integers in $[a,b]$.

\begin{Lemma}\label{g-regular 4g}
There exists a $g$-regular equi-difference $1$-D $(4g,3,2,1)$-OOC with $\lceil g/2\rceil$ codewords for any positive integer $g$, whose difference leave is $[1,g-1]_o\cup [3g+1,4g-1]_o\cup (4\cdot[1,g-1])$.
\end{Lemma}

\proof Let $${\cal C}=\left\{
\begin{array}{ll}
\{\{0,i,2i\}:i\in[g+1,2g-1]_o\}, &  {\rm if\ } g\equiv 0\ ({\rm mod}\ 2);\\
\{\{0,i,2i\}:i\in[g,2g-1]_o\}, &  {\rm if\ } g\equiv 1\ ({\rm mod}\ 2).\\
\end{array}
\right.
$$
Note that $\{0,4g/3,8g/3\}\not\in{\cal C}$, i.e., $i\neq 4g/3$ since $i$ is odd. Then ${\cal C}$ forms a $g$-regular equi-difference $1$-D $(4g,3,2,1)$-OOC with $\lceil g/2\rceil$ codewords, whose difference leave is $[1,g-1]_o\cup [3g+1,4g-1]_o\cup (4\cdot[1,g-1])$. \qed

\begin{Lemma}\label{cor:4 times}
Let $m\equiv 0\ ({\rm mod}\ 4)$. If there exists an optimal equi-difference $1$-D $(m/4,3,2$, $1)$-OOC with $\Psi^e(m/4,3,2,1)$ codewords, whose difference leave is $L$, then there exists an optimal equi-difference $1$-D $(m,3,2,1)$-OOC with $\lceil m/8\rceil+\Psi^e(m/4,3,2,1)$ codewords, whose difference leave is $(4\cdot L)\cup [1,m/4-1]_o\cup [3m/4+1,m-1]_o$.
\end{Lemma}

\proof By Lemma \ref{g-regular 4g}, there exists an $(m/4)$-regular equi-difference $1$-D $(m,3,2,1)$-OOC with $\lceil m/8\rceil$ codewords for any $m\equiv 0\ ({\rm mod}\ 4)$, whose difference leave is $ [1,m/4-1]_o\cup [3m/4+1,m-1]_o\cup (4\cdot[1,m/4-1])$. Then apply Construction \ref{filling g-regular} with an optimal equi-difference $1$-D $(m/4,3,2,1)$-OOC with $\Psi^e(m/4,3,2,1)$ codewords, whose difference leave is $L$, to obtain an equi-difference $1$-D $(m,3,2,1)$-OOC with $\lceil m/8\rceil+\Psi^e(m/4,3,2,1)$ codewords, whose difference leave is $(4\cdot L)\cup [1,m/4-1]_o\cup [3m/4+1,m-1]_o$. The optimality is ensured by Lemma \ref{bounf equi-diff}. \qed

\begin{Theorem}\label{lem:app 4 times-1}
There exists an optimal equi-difference $1$-D $(4^s r,3,2,1)$-OOC with
$$\Psi^e(4^s r,3,2,1)=(2^{2s+1}r+r-6)/12$$
codewords for any $s\geq 0$ and $r\equiv 2\ ({\rm mod}\ 4)$, whose difference leave is
$$\{2^{2s-1}r\}\cup \left(\bigcup_{i=1}^s (4^{s-i}\cdot ([1,4^{i-1}r-1]_o\cup[4^{i-1}3r+1,4^ir-1]_o)) \right).$$
\end{Theorem}

\proof We use induction on $s$. When $s=0$, the conclusion follows from Lemma \ref{lem:equi-2 mod 4}. Assume that the conclusion holds for $s=k-1$ ($k\geq 1$), i.e., there exists an optimal equi-difference $1$-D $(4^{k-1} r,3,2,1)$-OOC with $(2^{2(k-1)+1}r+r-6)/12$ codewords. Then apply Lemma \ref{cor:4 times} to obtain an optimal equi-difference $1$-D $(4^{k} r,3,2,1)$-OOC with $\lceil 4^{k-1} r/2\rceil+(2^{2(k-1)+1}r+r-6)/12=(2^{2k+1}r+r-6)/12$ codewords. One can check the difference leave for any given $s$ by induction. \qed

The difference leave of each optimal equi-difference $1$-D $(4^s r,3,2,1)$-OOC constructed in Theorem \ref{lem:app 4 times-1} contains $2^{2s-1}r$, which is a half of $4^s r$. In Section 5.2 we shall present  direct constructions for optimal $2$-D $(3\times m,3,2,1)$-OOCs which must contain three optimal equi-difference $1$-D $(m,3,2,1)$-OOCs as subcodes. However, we hope the three equi-difference $1$-D OOCs can use up the three half differences in their codewords. For this purpose, we present the following theorem to show optimal equi-difference $1$-D $(4^s r,3,2,1)$-OOCs whose difference leave do not contain $2^{2s-1}r$.

\begin{Theorem}\label{lem:app 4 times-11}
There exists an optimal equi-difference $1$-D $(4^s r,3,2,1)$-OOC with
$$\Psi^e(4^s r,3,2,1)=(2^{2s+1}r+r-6)/12$$
codewords for any $s\geq 1$ and $r\equiv 2\ ({\rm mod}\ 4)$, whose difference leave is
$$\{2^{2s-3}3r,2^{2s-3}5r\}\cup \left(\bigcup_{i=1}^s (4^{s-i}\cdot ([1,4^{i-1}r-1]_o\cup[4^{i-1}3r+1,4^ir-1]_o)) \right).$$
\end{Theorem}

\proof When $s=1$, by the proof of Theorem \ref{lem:app 4 times-1}, we can list all $(3r-2)/4$ codewords of an optimal equi-difference $1$-D $(4r,3,2,1)$-OOC as follows:
$${\cal A}=\{\{0,i,2i\}:i\in[r+1,2r-1]_o\}\cup (4\cdot \{\{0,i,2i\}:i=1,3,\ldots,r/2-2\}),$$
whose difference leave is $\{2r\}\cup [1,r-1]_o\cup[3r+1,4r-1]_o$. Let
$${\cal B}=({\cal A}\setminus \{\{0,3r/2,3r\}\})\cup \{\{0,r,2r\}\}.$$
Then ${\cal B}$ is an optimal equi-difference $1$-D $(4r,3,2,1)$-OOC, whose difference leave is $\{3r/2,5r/2\}\cup [1,r-1]_o\cup[3r+1,4r-1]_o$. Start from $\cal B$ and use induction on $s$. Then we can complete the proof by similar argument to that in Theorem \ref{lem:app 4 times-1}. \qed

\subsection{Constructions from conflict-avoiding codes}

Recall that a $1$-D $(m,k,k,1)$-OOC is called a conflict-avoiding code, denoted by a CAC$(m,k)$.

A CAC$(m,3)$ is said to be \emph{equi-difference} if each of its codewords is of the form $X=\{0,a,2a\}$, i.e., $|supp(\Delta X)|=2,3$ or $4$. Let $M^e(m,3)$ denote the largest size of codes among all equi-difference CAC$(m,3)$s for given $m$. An equi-difference CAC$(m,3)$ is said to be \emph{optimal} if it contains $M^e(m,3)$ codewords.

Clearly, when $m\not\equiv 0\ ({\rm mod}\ 3)$, an optimal equi-difference CAC$(m,3)$ is also an optimal equi-difference $1$-D $(m,3,2,1)$-OOC. Thus many known constructions for optimal equi-difference CAC$(m,3)$ in \cite{fls,lmsj,mzs,m0,wf} can be applied to constructions for optimal equi-difference $1$-D $(m,3,2,1)$-OOC.

An equi-difference CAC$(m,3)$ $\mathcal{C}$ is said to be \emph{tight} if $\bigcup_{X\in\mathcal{C}}supp(\Delta X)=Z_m\setminus\{0\}$. A tight equi-difference CAC$(m,3)$ is optimal. Momihara \cite{m0} gave a necessary and sufficient condition for the existence of a tight equi-difference CAC$(m,3)$. In Fu et al.\cite{fls}, the condition is restated in terms of multiplicative order of $2$ modulo $p$ for all prime factors $p$ of $m$.

\begin{Lemma}{\rm \cite{fls,m0}}\label{equi-cac1}
There exists a tight equi-difference CAC$(m,3)$ if and only if $m=4$ or $m\geq 3$ and $m=3^f m_0$ for $f\in\{0,1\}$, where any prime factor $p$ of $m_0$ satisfies $p\equiv 1\ ({\rm mod}\ 4)$ and $4|ord_{p}(2)$ whenever $p\equiv 1\ ({\rm mod}\ 8)$. Furthermore, a tight equi-difference CAC$(m,3)$ contains $1$ codeword for $m=4$, $(m-1)/4$ codewords for admissible $m\equiv 1,5\ ({\rm mod}\ 12)$, and $(m+1)/4$ codewords for admissible $m\equiv 3\ ({\rm mod}\ 12)$.
\end{Lemma}

We remark that the conditions on $m_0$ in Lemma \ref{equi-cac1} are fairly complex and one has to examine each prime factor of $m_0$. For this reason, only a few explicit series of odd $m$ are known (see \cite{lmsj,wf}).

\begin{Theorem}\label{lem:tight-1}
Let $r\equiv 1,5\ ({\rm mod}\ 12)$ satisfying that for any prime factor $p$ of $r$, $p\equiv 1\ ({\rm mod}\ 4)$ and $4|ord_{p}(2)$ whenever $p\equiv 1\ ({\rm mod}\ 8)$. Then
\begin{enumerate}
\item[$(1)$] there is an optimal equi-difference $1$-D $(r,3,2,1)$-OOC with $$\Psi^e(m,3,2,1)=(r-1)/4$$ codewords, whose difference leave is empty;
\item[$(2)$] there is an optimal equi-difference $1$-D $(4^s r,3,2,1)$-OOC with $$\Psi^e(4^s r,3,2,1)=(2^{2s-1}-2)r/3+(3r+1)/4$$ codewords for any $s\geq 1$, whose difference leave is
    $$\bigcup_{i=1}^s (4^{s-i}\cdot ([1,4^{i-1}r-1]_o\cup[4^{i-1}3r+1,4^ir-1]_o)).$$
\end{enumerate}
Note that it is allowed that $r=1$.
\end{Theorem}

\proof Since $r\not\equiv 0\ ({\rm mod}\ 3)$, (1) is straightforward by Lemmas \ref{bounf equi-diff} and  \ref{equi-cac1}. To prove (2), we use induction on $s$. When $s=1$, take an optimal equi-difference $1$-D $(r,3,2,1)$-OOC with $(r-1)/4$ codewords from (1), whose difference leave is empty. Then apply Lemma \ref{cor:4 times} to obtain an optimal equi-difference $1$-D $(4r,3,2,1)$-OOC with $(3r+1)/4$ codewords, whose difference leave is $[1,r-1]_o\cup [3r+1,4r-1]_o$. Assume that the conclusion holds for $s=k-1$ ($k\geq 2$). Then apply Lemma \ref{cor:4 times} again to complete the proof. \qed

\begin{Theorem}\label{lem:tight-2}
Let $r\equiv 3\ ({\rm mod}\ 12)$. If for any prime factor $p$ of $r/3$, $p\equiv 1\ ({\rm mod}\ 4)$ and $4|ord_{p}(2)$ whenever $p\equiv 1\ ({\rm mod}\ 8)$, then
\begin{enumerate}
\item[$(1)$] there is an optimal equi-difference $1$-D $(r,3,2,1)$-OOC with $$\Psi^e(m,3,2,1)=(r-3)/4$$ codewords, whose difference leave is $\{r/3,2r/3\}$;
\item[$(2)$] there is an optimal equi-difference $1$-D $(4^s r,3,2,1)$-OOC with $$\Psi^e(4^s r,3,2,1)=(2^{2s-1}-2)r/3+(3r-1)/4$$ codewords  for any $s\geq 1$, whose difference leave is $$\{2^{2s}r/3,2^{2s+1}r/3\}\cup \left(\bigcup_{i=1}^s (4^{s-i}\cdot ([1,4^{i-1}r-1]_o\cup[4^{i-1}3r+1,4^ir-1]_o)) \right).$$
\end{enumerate}
Note that it is allowed that $r=3$.
\end{Theorem}

\proof (1) By Lemma \ref{equi-cac1}, a tight equi-difference CAC$(r,3)$ with $(r+1)/4$ codewords exists for any $r\equiv 3\ ({\rm mod}\ 12)$ satisfying the assumption. Since $r$ is odd, such a tight CAC must contain the codeword $\{0,r/3,2r/3\}$. It follows that all codewords of the CAC except for the codeword $\{0,r/3,2r/3\}$ constitute an equi-difference $1$-D $(r,3,2,1)$-OOC with $(r-3)/4$ codewords, which is optimal by Lemma \ref{bounf equi-diff}.

(2) We use induction on $s$. When $s=1$, take an optimal equi-difference $1$-D $(r,3,2,1)$-OOC with $(r-3)/4$ codewords from (1), whose difference leave is $\{r/3,2r/3\}$. Then apply Lemma \ref{cor:4 times} to obtain an optimal equi-difference $1$-D $(4r,3,2,1)$-OOC with $(3r-1)/4$ codewords, whose difference leave is $\{4r/3,8r/3\}\cup [1,r-1]_o\cup [3r+1,4r-1]_o$. Assume that the conclusion holds for $s=k-1$ ($k\geq 2$). Then apply Lemma \ref{cor:4 times} again to complete the proof. \qed

\begin{Lemma}{\rm \cite{mzs}}\label{equi-cac-Ma}
Let $p\geq 5$ be any prime. There exists an optimal equi-difference CAC$(p,3)$ with
$$M^e(p,3)=\frac{p-1}{2^{ord_p(2)\ ({\rm mod}\ 2)}ord_p(2)}\times \lfloor \frac{1}{2} \times \frac{ord_p(2)}{2^{(ord_p(2)+1) \ ({\rm mod}\ 2)}} \rfloor $$
codewords.
\end{Lemma}

\begin{Theorem}\label{lem:non-tight-from Ma}
Let $p\geq 5$ be any prime and $M^e(p,3)$ be as in Lemma $\ref{equi-cac-Ma}$. Then
\begin{enumerate}
\item[$(1)$] there is an optimal equi-difference $1$-D $(p,3,2,1)$-OOC with $$\Psi^e(p,3,2,1)=M^e(p,3)$$ codewords;
\item[$(2)$] there is an optimal equi-difference $1$-D $(4^s p,3,2,1)$-OOC with $$\Psi^e(4^s p,3,2,1)=(2^{2s-1}-2)p/3+(p+1)/2+M^e(p,3)$$ codewords for any $s\geq 1$.
\end{enumerate}
\end{Theorem}

\proof Since $p\not\equiv 0\ ({\rm mod}\ 3)$, (1) follows immediately from Lemma \ref{equi-cac-Ma}. To prove (2), one can use induction on $s$ and apply Lemma \ref{cor:4 times} repeatedly. \qed

\section{Determination of $\Phi(n\times m,3,2,1)$ with $m\equiv 0\pmod{4}$}

\subsection{The cases of $n=1$ and $2$}

\begin{Lemma}\label{bound}
There exists an optimal $1$-D $(m,3,2,1)$-OOC for any $m\equiv 0\ ({\rm mod}\ 4)$ with
$$\Phi(1\times m,3,2,1)=\left\{
\begin{array}{lll}
\lfloor 7m/32\rfloor, & {\rm if\ } m\equiv 0\ ({\rm mod}\ 8),\\
\lfloor (7m+4)/32\rfloor, & {\rm if\ } m\equiv 4\ ({\rm mod}\ 8),\\
\end{array}
\right.
$$

\noindent codewords with the exception of $\Phi(1\times 64,3,2,1)=13.$

\end{Lemma}

\proof For $m=48$, we give an explicit construction for a $1$-D $(48,3,2,1)$-OOC with $10$ codewords as follows, which is defined on $Z_{48}$. Lemma \ref{bound1} ensures its optimality.
\begin{center}
\begin{tabular}{llllll}
 $\{0,3,6\},$ &$\{0,7,14\},$& $\{0,11,22\},$& $\{0,15,30\},$ &$\{0,19,38\},$ &
 $\{0,23,46\},$\\ $\{0,1,17\},$& $\{0,5,9\},$ &$\{0,13,21\},$ &$\{0,12,24\}$.&&
\end{tabular}
\end{center}
For $m\equiv 0\ ({\rm mod}\ 4)$ and $m\neq 48$, let $m=4u$ and $u\neq 12$. Write $T=\{0,m/3,2m/3\}$ when $m\equiv 0\ ({\rm mod}\ 3)$. Then an optimal $1$-D $(m,3,2,1)$-OOC with $\Phi(1\times m,3,2,1)$ codewords is constructed in the following table. Note that when $m\not\equiv 0\ ({\rm mod}\ 3)$, an optimal $1$-D $(m,3,3,1)$-OOC (or CAC$(m,3)$) is also an optimal $1$-D $(m,3,2,1)$-OOC.
\begin{center}
\begin{tabular}{|c|c|c|}\hline
\multicolumn{2}{|c|}{$u$} & Source \\\hline
\multirow{3}{*}{$0\ ({\rm mod}\ 8)$} & $0,8\ ({\rm mod}\ 32)$ &  Construction 3.1 in \cite{mfu} \\\cline{2-3}
& $16,24\ ({\rm mod}\ 32)$, $\neq16$ &  Construction 3.2 in \cite{mfu}\\\cline{2-3}
& $16$ &  Theorem \ref{cac} \\\hline

\multirow{2}{*}{$4\ ({\rm mod}\ 8)$, $\neq12$} & $4,20\ ({\rm mod}\ 24)$ &  Theorem \ref{cac}\\\cline{2-3}
& $12\ ({\rm mod}\ 24)$, $\neq12$ & Constructions 3.6, 3.7 in \cite{mfu}; discard $T$ \\\hline

\multirow{2}{*}{$2\ ({\rm mod}\ 8)$} & $2,10\ ({\rm mod}\ 24)$ &  Theorem \ref{cac}\\\cline{2-3}
& $18\ ({\rm mod}\ 24)$ & Constructions 5.3, 5.4 in \cite{mfu}; discard $T$ \\\hline

\multirow{2}{*}{$6\ ({\rm mod}\ 8)$} & $22,30\ ({\rm mod}\ 32)$ &  Construction 5.5 in \cite{mfu}\\\cline{2-3}
& $6,14\ ({\rm mod}\ 32)$ & Construction 5.6 in \cite{mfu} \\\hline

$1\ ({\rm mod}\ 8)$  &  & Construction 3.1, 3.2 in \cite{fu}\\\hline

\multirow{2}{*}{$3\ ({\rm mod}\ 8)$} & $3\ ({\rm mod}\ 24)$ &  Construction 3.3 in \cite{fu}; discard $T$ \\\cline{2-3}
&  $11,19\ ({\rm mod}\ 24)$ & Theorem \ref{cac} \\\hline

\multirow{2}{*}{$5\ ({\rm mod}\ 8)$} & $21\ ({\rm mod}\ 24)$ &  Construction 3.7, 3.8, 3.9 in \cite{fu}; discard $T$ \\\cline{2-3}
&  $5,13\ ({\rm mod}\ 24)$ & Theorem \ref{cac} \\\hline

$7\ ({\rm mod}\ 8)$  &  & Construction 3.10 in \cite{fu}\\ \hline
\end{tabular}
\end{center}

It should be noticed that, when $u\equiv0,8 \pmod{32}$, although Construction $3.1$ in \cite{mfu}
was only used for the cases of $u\equiv 8,32,40,64\pmod{96}$, it is readily checked that the same codewords listed in Construction $3.1$ in \cite{mfu} can also produce our required optimal OOCs for $u\equiv0,72 \pmod{96}$. The similar things happen when $u\equiv16,24\pmod{32}$ and $u\equiv6\pmod{8}$.

We give another two examples to illustrate how to use the table. When $u\equiv 3\ ({\rm mod}\ 24)$, in Construction $3.3$ of \cite{fu}, an optimal CAC$(m,3)$ with $(7m+12)/32$ codewords is constructed, where $m=4u$ and $T=\{0,m/3,2m/3\}$ is one of codewords. Then all codewords of the CAC$(m,3)$ except for the codeword $T$ constitute an optimal $1$-D $(m,3,2,1)$-OOC with $(7m-20)/32$ codewords.

When $u\equiv 7\ ({\rm mod}\ 8)$, in Construction $3.10$ of \cite{fu}, an optimal CAC$(m,3)$ with $ (7m-4)/32$ codewords is constructed, where $m=4u$ and $T=\{0,m/3,2m/3\}$ is not a codeword. Thus this CAC$(m,3)$ is also an optimal $1$-D $(m,3,2,1)$-OOC with $(7m-4)/32$ codewords. \qed

\begin{Lemma}\label{n=2}
There exists an optimal $2$-D $(2\times m,3,2,1)$-OOC with $\Phi(2\times m,3,2,1)=3m/4$ codewords for any $m\equiv 0\ ({\rm mod}\ 4)$, with the exception of $\Phi(2\times 4,3,2,1)=2$.
\end{Lemma}

\proof We construct the required codes on $I_2\times Z_m$. When $m\equiv 0\ ({\rm mod}\ 8)$ and $m\geq8$, the required $3m/4$ codewords are
\begin{center}
\begin{tabular}{ll}
$\{(0,0),(0,i),(0,2i)\}$, &$i\in[3,m/4-1]_o\cup\{m/2-1\}$, ($i=3$ if $m=8$);\\
$\{(1,0),(1,i),(1,2i)\}$, &$i\in[m/4+1,m/2-1]_o$;\\
$\{(0,0),(0,1+2i),(1,m/4-1+i)\}$,& $i\in[m/8,m/4-2]$, (null if $m=8$);\\
$\{(1,0),(1,1+2i),(0,3m/4+2+i)\}$,& $i\in[0,m/8-1]$;\\
$\{(0,0),(0,4+4i),(1,3m/4+1+2i)\}$,& $i\in[0,m/8-1]$;\\
$\{(1,0),(1,4+4i),(0,m/4+4+2i)\}$, &$i\in[0,m/8-1]$;\\
$\{(0,0),(0,1),(1,3m/4-1)\}$.&\\
\end{tabular}
\end{center}
When $m\equiv 4\ ({\rm mod}\ 8)$ and $m\geq12$, the required $3m/4$ codewords are
\begin{center}
\begin{tabular}{ll}
$\{(0,0),(0,i),(0,2i)\}$, &$i\in[m/4+2,m/2-3]_o$, (null if $m=12$);\\
$\{(1,0),(1,i),(1,2i)\}$,& $i\in[1,m/4]_o$;\\
$\{(0,0),(0,1+2i),(1,m/4+i)\}$, &$i\in[0,(m-4)/8]\setminus\{1\}$;\\
$\{(1,0),(1,1+2i),(0,3m/4+1+i)\}$,& $i\in[(m+4)/8,m/4-1]$;\\
$\{(0,0),(0,4+4i),(1,3m/4+1+2i)\}$,& $i\in[1,(m-12)/8]$, (null if $m=12$);\\
$\{(1,0),(1,4+4i),(0,m/4+2+2i)\}$,& $i\in[0,(m-12)/8]$;\\
$\{(0,0),(0,m/2-1),(1,m/4-2)\}$,& $\{(0,0),(0,3),(1,3m/4)\}$,\\
$\{(0,0),(0,m/2),(1,3m/4+1)\}$, &$\{(0,0),(0,2),(0,4)\}$.\\
\end{tabular}
\end{center}
When $m=4$, the required two codewords are $\{(0,0),(0,1),(0,2)\}$ and  $\{(1,0),(1,1),(1,2)\}$. Lemmas \ref{n=2bound} and \ref{2*4} ensure the optimality of these codes. \qed

\subsection{The case of $n=3$} \label{sec:n=3}

This section is devoted to constructing optimal $2$-D $(3\times m,3,2,1)$-OOCs. By Lemma \ref{0mod4}, $\Phi(3\times m,3,2,1)\leq 3m/2+\Psi^e(m,3,2,1)$, and so by Remark \ref{remark}, any optimal $2$-D $(3\times m,3,2,1)$-OOC must contain $3$ optimal equi-difference $1$-D $(m,3,2,1)$-OOCs as subcodes.

\subsubsection{$m\equiv 8\ ({\rm mod}\ 16)$}

\begin{Lemma}\label{3*8}
There exists an optimal $2$-D $(3\times 8,3,2,1)$-OOC with $\Phi(3\times 8,3,2,1)=13$ codewords.
\end{Lemma}

\proof The required OOC is constructed on $I_3\times Z_8$ as follows:
\begin{center}\tabcolsep 0.02in
\begin{tabular}{llll}
$\{(0,0),(0,2),(0,4)\}$,&
$\{(1,0),(1,2),(1,4)\}$,&
$\{(2,0),(2,2),(2,4)\}$;\\
$\{(0,0),(0,1),(1,6)\}$,&
$\{(0,0),(0,3),(1,7)\}$,&
$\{(1,0),(1,1),(2,5)\}$,&
$\{(1,0),(1,3),(2,3)\}$,\\
$\{(0,0),(2,5),(2,6)\}$,&
$\{(0,0),(2,4),(2,7)\}$;\\
$\{(0,0),(1,2),(2,0)\}$,&
$\{(0,0),(1,3),(2,2)\}$,&
$\{(0,0),(1,0),(2,1)\}$,&
$\{(0,0),(1,1),(2,3)\}$.
\end{tabular}
\end{center}
The optimality is ensured by Lemma \ref{0mod4} and Theorem \ref{lem:app 4 times-11}. Note that each of the first three codewords can be seen as an optimal equi-difference $1$-D $(8,3,2,1)$-OOC, which is defined on $\{x\}\times Z_8$ for some $x\in\{0,1,2\}$. The middle six codewords used up all the remaining pure differences which are not from the first three codewords. All mixed differences are also used up. \qed

\begin{Remark}\label{remark:difficult}
Lemma $\ref{3*8}$ helps us to understand the structure of codewords of optimal $2$-D $(3\times m,3,2,1)$-OOCs with $3m/2+\Psi^e(m,3,2,1)$ codewords $($if it exists$)$. On one hand, such kind of OOCs must contain three optimal equi-difference $1$-D $(m,3,2,1)$-OOCs as subcodes. On the other hand, all pure differences and mixed differences must be used up. The two facts make it difficult to find effective recursive constructions, especially filling constructions, for optimal $2$-D $(3\times m,3,2,1)$-OOCs.
\end{Remark}

In the following we present three infinite families of optimal $2$-D $(3\times m,3,2,1)$-OOCs via direct constructions.

\begin{Lemma}\label{8 mod 16}
There exists an optimal $2$-D $(3\times m,3,2,1)$-OOC with $\Phi(3\times m,3,2,1)=(27m-8)/16$ codewords  for any $m\equiv 8\ ({\rm mod}\ 16)$.
\end{Lemma}

\proof For $m\equiv 8\ ({\rm mod}\ 16)$, by Lemma \ref{0mod4}, $\Phi(3\times m,3,2,1)\leq 3m/2+\Psi^e(m,3,2,1)$, and by Theorem \ref{lem:app 4 times-1}, $\Psi^e(m,3,2,1)=(3m-8)/16$. So $\Phi(3\times m,3,2,1)\leq (27m-8)/16$.

When $m\equiv 8\ ({\rm mod}\ 16)$ and $m\geq 24$, the required $(27m-8)/16$ codewords are divided into two parts. The first part consists of $(9m-24)/16$ codewords:
$$\{(x,0),(x,a),(x,2a)\},\ \ x\in \{0,1,2\}\ {\rm and}\ \{0,a,2a\}\in \mathcal{B},$$
where $\mathcal{B}$ is an optimal equi-difference $1$-D $(m,3,2,1)$-OOC with $(3m-8)/16$ codewords, whose difference leave is $\{3m/8,5m/8\}\cup [1,m/4-1]_o\cup [3m/4+1,m-1]_o$ (see Theorem \ref{lem:app 4 times-11} by taking $s=1$ and $r=m/4$). The second part consists of $9m/8+1$ codewords:
\begin{center}
\begin{tabular}{lll}
$\{(0,0),(0,1+2i),(1,7m/8+i)\}$, & $i\in[0,m/8-1]$; \\
$\{(1,0),(1,1+2i),(2,m/2+2+i)\}$,  & $i\in[0,m/8-1]$; \\
$\{(0,0),(2,7m/8-3-i),(2,7m/8+i)\}$, & $i\in[0,m/8-2]$; \\
$\{(0,0),(0,3m/8),(1,3m/4-1)\}$, & $\{(1,0),(1,3m/8),(2,3m/4)\}$, \\
$\{(0,0),(2,m-1),(2,0)\}$,  &  $\{(0,0),(2,7m/8-2),(2,m/4-2)\}$;\\
$\{(0,0),(1,i),(2,1+2i)\}$, & $i\in[0,3m/8-2]$; \\
$\{(0,0),(1,3m/8+i),(2,2+2i)\}$, & $i\in[0,3m/8-2]\setminus\{m/8-2\}$; \\
$\{(0,0),(1,m/2-2),(2,7m/8-1)\}$.
\end{tabular}
\end{center}
When $m=8$, the conclusion follows from Lemma \ref{3*8}. \qed

\subsubsection{$m\equiv 32\ ({\rm mod}\ 64)$}

\begin{Lemma}\label{32 mod 64}
There exists an optimal $2$-D $(3\times m,3,2,1)$-OOC with $\Phi(3\times m,3,2,1)=(107m-32)/64$ codewords for any $m\equiv 32\ ({\rm mod}\ 64)$.
\end{Lemma}

\proof For $m\equiv 32\ ({\rm mod}\ 64)$, by Lemma \ref{0mod4}, $\Phi(3\times m,3,2,1)\leq 3m/2+\Psi^e(m,3,2,1)$, and by Theorem \ref{lem:app 4 times-1}, $\Psi^e(m,3,2,1)=(11m-32)/64$. So $\Phi(3\times m,3,2,1)\leq (107m-32)/64$.

When $m=32$, an optimal $2$-D $(3\times 32,3,2,1)$-OOC with $53$ codewords is listed as follows:
\begin{center}
\hspace*{-2cm}\begin{tabular}{lllllll}
$\{(x,0),(x,a),(x,2a)\}$,\ \ $x\in\{0,1,2\}$ and $a\in\{8,9,11,13,15\}$;\\
\end{tabular}
\begin{tabular}{lllllll}
$\{(0,0),     (0,12),     (1,25)\}$,&
$\{(0,0),     (0,1),     (1,6)\}$,&
$\{(0,0),     (0,3),     (1,7)\}$,\\
$\{(0,0),     (0,5),     (1,8)\}$,&
$\{(0,0),     (0,7),     (1,9)\}$,&
$\{(0,0),     (0,4),     (1,31)\}$,\\

$\{(1,0),     (1,1),     (2,14)\}$,&
$\{(1,0),     (1,12),     (2,23)\}$,&
$\{(1,0),     (1,3),     (2,20)\}$,\\
$\{(1,0),     (1,5),     (2,21)\}$,&
$\{(1,0),     (1,7),     (2,22)\}$,&
$\{(1,0),     (1,4),     (2,31)\}$,\\

$\{(0,0),     (2,9),     (2,21)\}$,&
$\{(0,0),     (2,0),     (2,1)\}$,&
$\{(0,0),     (2,7),     (2,10)\}$,\\
$\{(0,0),     (2,6),     (2,11)\}$,&
$\{(0,0),     (2,5),     (2,12)\}$,&
$\{(0,0),     (2,22),     (2,26)\}$;\\

$\{(0,0),     (1,0),     (2,8)\}$,&
$\{(0,0),     (1,1),     (2,4)\}$,&
$\{(0,0),     (1,10),     (2,28)\}$,\\
$\{(0,0),     (1,19),     (2,31)\}$,&
$\{(0,0),     (1,21),     (2,30)\}$,&
$\{(0,0),     (1,22),     (2,29)\}$,\\
$\{(0,0),     (1,12),     (2,14)\}$,&
$\{(0,0),     (1,14),     (2,20)\}$,&
$\{(0,0),     (1,15),     (2,19)\}$,\\
$\{(0,0),     (1,16),     (2,16)\}$,&
$\{(0,0),     (1,17),     (2,18)\}$,&
$\{(0,0),     (1,18),     (2,23)\}$,\\
$\{(0,0),     (1,11),     (2,3)\}$,&
$\{(0,0),     (1,20),     (2,13)\}$,&
$\{(0,0),     (1,23),     (2,17)\}$,\\
$\{(0,0),     (1,24),     (2,2)\}$,&
$\{(0,0),     (1,26),     (2,24)\}$,&
$\{(0,0),     (1,28),     (2,15)\}$,\\
$\{(0,0),     (1,29),     (2,25)\}$,&
$\{(0,0),     (1,30),     (2,27)\}$.
\end{tabular}
\end{center}
Note that for any $x\in \{0,1,2\}$, $\{\{(x,0),(x,a),(x,2a)\}:a\in\{8,9,11,13,15\}\}$
forms an optimal equi-difference $1$-D $(32,3,2,1)$-OOC defined on $\{x\}\times Z_{32}$, whose difference leave is $\{1,3,4,5,7,12,20,25,27,28,29,31\}$ (see Theorem \ref{lem:app 4 times-11} by taking $s=2$ and $r=2$).

When $m\equiv 32\ ({\rm mod}\ 64)$ and $m\geq 96$, the required $(107m-32)/64$ codewords are divided into two parts. The first part consists of $(33m-96)/64$ codewords:
$$\{(x,0),(x,a),(x,2a)\},\ \ x\in \{0,1,2\}\ {\rm and}\ \{0,a,2a\}\in \mathcal{B},$$
where $\mathcal{B}$ is an optimal equi-difference $1$-D $(m,3,2,1)$-OOC with $(11m-32)/64$ codewords, whose difference leave is $\{3m/8,5m/8\}\cup [1,m/4-1]_o\cup [3m/4+1,m-1]_o\cup (4\cdot ([1,m/16-1]_o\cup[3m/16+1,m/4-1]_o))$ (see Theorem \ref{lem:app 4 times-11} by taking $s=2$ and $r=m/16$). The second part consists of $37m/32+1$ codewords:
\begin{center}\tabcolsep 0.015in
\hspace*{-3.9cm}\begin{tabular}{ll}
$\{(0,0),(0,1+2i),(1,m/8+2+i)\}$,& $i\in[0,m/8-1]$;\\
$\{(0,0),(0,4+8i),(1,7m/8+3+4i)\}$, &$i\in[0,m/32-1]$;\\
$\{(1,0),(1,3+2i),(2,5m/8+i)\}$, & $i\in[0,m/8-2]$;\\
$\{(1,0),(1,4+8i),(2,7m/8+3+4i)\}$, & $i\in[0,m/32-1]$; \\
$\{(0,0),(2,m/4-1-i),(2,m/4+2+i)\}$, & $i\in[0,m/8-2]$; \\
$\{(0,0),(2,3m/4-2-4i),(2,3m/4+2+4i)\}$, & $i\in[0,m/32-1]$; \\
\end{tabular}
\hspace*{-2.3cm}\begin{tabular}{lll}
$\{(0,0),(0,3m/8),(1,13m/16-1)\}$,\ &  $\{(1,0),(1,1),(2,7m/16)\}$,\\
$\{(1,0),(1,3m/8),(2,3m/4-1)\}$,\ &  $\{(0,0),(2,m/4+1),(2,5m/8+1)\}$,\\
$\{(0,0),(2,m/16-2),(2,m/16-1)\}$;
\end{tabular}
\begin{tabular}{ll}
$\{(0,0),(1,3m/4+2i),(2,3m/4-3-2i)\}$, & $i\in[0,m/16-1]\setminus\{m/16-2\}$; \\
$\{(0,0),(1,7m/8+2i),(2,5m/8+4i)\}$, & $i\in[0,m/16-1]$; \\
$\{(0,0),(1,3m/4+1+4i),(2,3m/4-1+2i)\}$, & $i\in[0,m/16-1]\setminus\{(m-32)/64\}$; \\
$\{(0,0),(1,m/4+2+2i),(2,3m/8+1+i)\}$, & $i\in[0,m/8-2]$; \\
$\{(0,0),(1,m/4+3+2i),(2,m/2+1+i)\}$, & $i\in[0,m/8-3]\setminus\{3m/32-2\}$; \\
$\{(0,0),(1,m/2+4+2i),(2,1+i)\}$, & $i\in[0,m/8-3]\setminus\{m/16-3,m/16-2\}$; \\
$\{(0,0),(1,m/2+1+2i),(2,7m/8-1+i)\}$, & $i\in[0,m/8-4]$; \\
\end{tabular}
\hspace*{-2cm}\begin{tabular}{lll}
$\{(0,0),(1,0),(2,0)\}$,\ & $\{(0,0),(1,1),(2,m/4)\}$,\\
$\{(0,0),(1,m/2-1),(2,m-3)\}$,\ & $\{(0,0),(1,m/2),(2,m/8-1)\}$,\\
$\{(0,0),(1,m/2+2),(2,m/8)\}$,\ & $\{(0,0),(1,3m/4-5),(2,m/2)\}$,\\
$\{(0,0),(1,3m/4-3),(2,m-2)\}$,\ & $\{(0,0),(1,3m/4-1),(2,m-1)\}$,\\
$\{(0,0),(1,7m/8-4),(2,m-4)\}$,\ & $\{(0,0),(1,5m/8-2),(2,25m/32-2)\}$,\\
$\{(0,0),(1,5m/8),(2,19m/32-1)\}$. \\
\end{tabular}
\end{center} 

\subsubsection{$m\equiv 4,20\ ({\rm mod}\ 48)$}

\begin{Lemma}\label{4}
There exists an optimal $2$-D $(3\times 4,3,2,1)$-OOC with $\Phi(3\times 4,3,2,1)=6$ codewords.
\end{Lemma}

\proof The required OOC is constructed on $I_3\times Z_4$ as follows:
\begin{center}
\begin{tabular}{llll}
$\{(x,0),(x,1),(x,2)\}$,&
$\{(0,0),(1,a),(2,b)\}$,
\end{tabular}
\end{center}
where $x\in \{0,1,2\}$ and $(a,b)\in\{(0,0),(1,3),(3,2)\}$. The optimality is ensured by Lemma \ref{3*4-bound}. \qed

\begin{Lemma}\label{4 mod 16}
Let $m>4$ and $m\equiv 4,20\ ({\rm mod}\ 48)$ satisfying that for any prime factor $p$ of $m/4$, $p\equiv 1\ ({\rm mod}\ 4)$ and $4|ord_{p}(2)$ whenever $p\equiv 1\ ({\rm mod}\ 8)$. Then there exists an optimal $2$-D $(3\times m,3,2,1)$-OOC with $\Phi(3\times m,3,2,1)=(27m+4)/16$ codewords.
\end{Lemma}

\proof For any $m$ given in the assumption, by Lemma \ref{0mod4}, $\Phi(3\times m,3,2,1)\leq 3m/2+\Psi^e(m,3,2,1)$, and by Theorem \ref{lem:tight-1}(2), $\Psi^e(m,3,2,1)=(3m+4)/16$. So $\Phi(3\times m,3,2,1)\leq (27m+4)/16$.

When $m=20$, an optimal $2$-D $(3\times 20,3,2,1)$-OOC with $34$ codewords is listed as follows:
\begin{center}
\hspace*{-2.95cm}\begin{tabular}{lllllll}
$\{(x,0),(x,i),(x,2i)\}$,\ \ $x\in \{0,1,2\}$ and $i\in\{4,5,7,9\}$;
\end{tabular}
\begin{tabular}{lllllll}
$\{(0,0),(0,1),(1,18)\}$,&
$\{(0,0),(0,3),(1,19)\}$,&
$\{(1,0),(1,1),(2,18)\}$,\\
$\{(1,0),(1,3),(2,19)\}$,&
$\{(0,0),(2,17),(2,18)\}$,&
$\{(0,0),(2,16),(2,19)\}$;\\
$\{(0,0),(1,a),(2,b)\}$,&
\end{tabular}
\end{center}
\noindent where $(a,b)\in\{(0,1),(1,3),(2,2),(3,11),(4,13),(5,10),(6,9),(7,14),(8,12),(9,15),(10$, $0),(11,4),(12,7),(13,5),(14,8),(15,6)\}$. Note that for any $x\in \{0,1,2\}$, $\{\{(x,0),(x,i)$, $(x,2i)\}:i\in\{4,5,7,9\}\}$
forms an optimal equi-difference $1$-D $(20,3,2,1)$-OOC defined on $\{x\}\times Z_{20}$, whose difference leave is $\{1,3,17,19\}$ (see Theorem \ref{lem:tight-1}(2) by taking $s=1$ and $r=5$).

When $m=52$, an optimal $2$-D $(3\times 52,3,2,1)$-OOC with $88$ codewords is listed as follows:
\begin{center}
\begin{tabular}{ll}
$\{(x,0),(x,i),(x,2i)\}$, & $x\in \{0,1,2\}$ and $i\in\{4,12,16,13,15,17,19,21,23,25\}$;\\
$\{(0,0),(0,1+2i),(1,46+i)\}$, & $i\in[0,5]$;\\
$\{(1,0),(1,1+2i),(2,46+i)\}$, & $i\in[0,5]$;\\
$\{(0,0),(2,45-i),(2,46+i)\}$, & $i\in[0,5]$;\\
$\{(0,0),(1,a),(2,b)\}$,
\end{tabular}
\end{center}
where $(a,b)\in\{(0,12), (2,6), (3,8), (4,14), (5,5), (6,7), (11,13), (1,16), (14,22), (16,19)$, $(18,24), (7,28), (12,25), (13,27), (22,29), (8,30), (9,32), (10,34), (20,31), (24,33), (15,35)$, $(17,36), (19,37), (21,38), (23,39), (33,15), (34,18), (27,0), (28,2), (29,4), (30,10), (32,11)$, $(25,1), (31,9), (26,3), (35,21), (36,17), (37,20), (38,23), (39,26)\}$. Note that for any $x\in \{0,1,2\}$, $\{\{(x,0),(x,i),(x,2i)\}:i\in\{4,12,16,13,15,17,19,21,23,25\}\}$
forms an optimal equi-difference $1$-D $(52,3,2,1)$-OOC defined on $\{x\}\times Z_{52}$, whose difference leave is $[1,11]_o\cup [41,51]_o$ (see Theorem \ref{lem:tight-1}(2) by taking $s=1$ and $r=13$).

When $m\equiv 4,20\ ({\rm mod}\ 48)$, $m\geq68$ and $m$ satisfies the condition in the assumption, the required $(27m+4)/16$ codewords are divided into three parts. The first part consists of $(9m+12)/16$ codewords:
$$\{(x,0),(x,a),(x,2a)\},\ \ x\in \{0,1,2\}\ {\rm and}\ \{0,a,2a\}\in \mathcal{B},$$
where $\mathcal{B}$ is an optimal equi-difference $1$-D $(m,3,2,1)$-OOC with $(3m+4)/16$ codewords, whose difference leave is $[1,m/4-1]_o\cup [3m/4+1,m-1]_o$ (see Theorem \ref{lem:tight-1}(2) by taking $s=1$ and $r=m/4$). The second part consists of $(3m+108)/8$ codewords:
\begin{center}
\begin{tabular}{lll}
$\{(0,0),(0,1+2i),(1,(7m+4)/8+i)\}$, & $i\in[0,(m-20)/8]$; \\
$\{(1,0),(1,1+2i),(2,m/2+i)\}$, & $i\in[0,(m-12)/8]$; \\
$\{(0,0),(2,(7m-12)/8-i),(2,(7m-4)/8+i)\}$, & $i\in[0,(m-12)/8]$; \\
\end{tabular}
\begin{tabular}{lll}
$\{(0,0),(0,m/4-2),(1,(11m-12)/16)\}$, & $\{(0,0),(1,(3m-4)/8),(2,(m-4)/16)\}$, \\
$\{(0,0),(1,(9m+12)/16),(2,m/2-1)\}$, &  $\{(0,0),(1,(3m+4)/8),(2,(3m+4)/16)\}$,\\
$\{(0,0),(1,3m/4+1),(2,2)\}$, &  $\{(0,0),(1,m-1),(2,3m/4-3)\}$,\\
$\{(0,0),(1,m/4),(2,m-1)\}$, &  $\{(0,0),(1,m/4-1),(2,m/4-3)\}$,\\
$\{(0,0),(1,3m/4),(2,0)\}$, &  $\{(0,0),(1,(3m+12)/8),(2,(m+12)/8)\}$,\\
$\{(0,0),(1,(5m-4)/8),(2,m/4-1)\}$, &  $\{(0,0),(1,(5m+4)/8),(2,m/4+1)\}$,\\
$\{(0,0),(1,3m/4-1),(2,(5m-20)/8)\}$, &  $\{(0,0),(1,m/2+1),(2,(3m+4)/8)\}$,\\
$\{(0,0),(1,m/2),(2,m/2)\}$, &  $\{(0,0),(1,m/2-1),(2,m/2-2)\}$.
\end{tabular}
\end{center}
The third part consists of $(3m-56)/4$ codewords. Let $T=[0,(3m-36)/8]\setminus\{(m-20)/16,(m-28)/8,(m-20)/8,(m-12)/8,(3m-28)/16,
m/4-3,m/4-2,(5m-52)/16\}$. If $m\equiv 4,68\ ({\rm mod}\ 96)$ and $m\geq 68$, then we take
\begin{center}
\begin{tabular}{lll}
$\{(0,0),(1,i),(2,1+2i)\}$, & $i\in[0,(3m-12)/8]\setminus\{(3m-12)/32,m/4-1,m/4\}$; \\
\end{tabular}
\begin{tabular}{lll}
$\{(0,0),(1,(3m-12)/32),(2,3m/4-1)\}$, & $\{(0,0),(1,(13m+12)/32),(2,m/2+1)\}$; \\
\end{tabular}
\begin{tabular}{lll}
$\{(0,0),(1,(3m+20)/8+i),(2,4+2i)\}$, & $i\in T\setminus \{(m-68)/32\}$.
\end{tabular}
\end{center}
If $m\equiv 20,52\ ({\rm mod}\ 96)$ and $m\geq 116$, then we take
\begin{center}
\begin{tabular}{lll}
$\{(0,0),(1,i),(2,1+2i)\}$, & $i\in[0,(3m-12)/8]\setminus\{(m-20)/32,m/4-1,m/4\}$; \\
\end{tabular}
\begin{tabular}{lll}
$\{(0,0),(1,(15m+20)/32),(2,m/2+1)\}$, & $\{(0,0),(1,(m-20)/32),(2,3m/4-1)\}$;
\end{tabular}
\begin{tabular}{lll}
$\{(0,0),(1,(3m+20)/8+i),(2,4+2i)\}$,  &  $i\in T\setminus\{(3m-60)/32\}$.
\end{tabular}
\end{center}

\subsection{A recursive construction from $m$-cyclic group divisible
designs}

Let $K$ be a set of positive integers. A {\em group divisible
design} (GDD) $K$-GDD is a triple ($X, {\cal G},{\cal A}$)
satisfying that ($1$) $\cal G$ is a partition
of a finite set $X$ into subsets (called {\em groups}); ($2$) $\cal
A$ is a set of subsets of $X$ (called {\em blocks}), each of
cardinality from $K$, such that every $2$-subset of $X$ is either
contained in exactly one block or in exactly one group, but not in
both. If $\cal G$ contains $u_i$ groups of size $g_i$ for $1\leq
i\leq r$, then we call $g_1^{u_1}g_2^{u_2}\cdots g_r^{u_r}$ the {\em type} of the GDD. If $K=\{k\}$, we write $k$-GDD instead of $\{k\}$-GDD.

An {\em automorphism group} of a GDD $(X,{\cal G},{\cal A})$ is a permutation group on $X$ leaving ${\cal G}$ and ${\cal A}$ invariant, respectively. Given an automorphism group of a GDD, all blocks of the GDD can be partitioned into some orbits under this automorphism group. Choose any fixed block from each orbit and call it a {\em base block} of the GDD.

Suppose $(X,\mathcal{G},\mathcal{A})$ is a $K$-GDD of type
$(v_{1}m)^{u_1} (v_{2}m)^{u_2}\cdots (v_{r}m)^{u_r}$. If its automorphism group contains a
permutation on $X$ that is the product of
$\sum_{i=1}^{r}v_iu_i$ disjoint $m$-cycles fixing each group of $\cal G$ and leaving $\mathcal{B}$ invariant, then this design is said to be \emph{$m$-cyclic}.

\begin{Lemma}{\rm\cite{wc2}} \label{h-cyclic gdd}
An $m$-cyclic $3$-GDD of type $(vm)^u$ exists if and only if $(1)$ when $u=3$, $m$ is odd, or $m$ is even and $v$ is even; $(2)$ when $u\geq 4$, $(u-1)vm\equiv0 \pmod{2}$,
$u(u-1)vm\equiv0 \pmod{3}$, and
$v\equiv0 \pmod{2}$ if
$u\equiv2,3 \pmod{4}$ and
$m\equiv2 \pmod{4}$.
\end{Lemma}

The following construction is a variation of Construction 4.6 in \cite{fwwz}.

\begin{Construction}\label{from 3-SCGDD} Suppose that there exist
\begin{enumerate}
\item[$(1)$] an $m$-cyclic $k$-GDD of type $(v_{1}m)^{u_1} (v_{2}m)^{u_2}\cdots (v_{r}m)^{u_r}$ with $b$ base blocks;
\item[$(2)$] a $2$-D $(v_i\times m,k,\lambda_a,1)$-OOC with $f_i$ codewords for each $1\leq i\leq r$.
\end{enumerate}
Then there exists a $2$-D $((\sum_{i=1}^r v_i u_i)\times m,k,\lambda_a,1)$-OOC with $b+\sum_{i=1}^r u_if_i$ codewords.
\end{Construction}

\begin{Lemma}\label{n-0-mod 3}
Let $m\equiv 0\ ({\rm mod}\ 4)$. If there is an optimal $2$-D $(3\times m,3,2,1)$-OOC with $3m/2+\Psi^e(m,3,2,1)$ codewords, then there is an optimal $2$-D $(n\times m,3,2,1)$-OOC with $n(nm+2\Psi^e(m,3,2,1))/6$ codewords for any $n\equiv 0\ ({\rm mod}\ 3)$ and $n\geq 12$.
\end{Lemma}

\proof By Lemma \ref{h-cyclic gdd}, there exists an $m$-cyclic $3$-GDD of type $(3m)^{n/3}$ for any $m\equiv 0\ ({\rm mod}\ 4)$, $n\equiv 0\ ({\rm mod}\ 3)$ and $n\geq 12$, which contains $mn(n-3)/6$ base blocks. Then apply Construction \ref{from 3-SCGDD} with an optimal $2$-D $(3\times m,3,2,1)$-OOC with $3m/2+\Psi^e(m,3,2,1)$ codewords to obtain a $2$-D $(n\times m,3,2,1)$-OOC with $n(nm+2\Psi^e(m,3,2,1))/6$ codewords, which is optimal by Lemma \ref{0mod4}. \qed

\begin{Corollary}\label{cor:=0 mod 3}
There is an optimal $2$-D $(n\times m,3,2,1)$-OOC with
$$\Phi(n\times m,3,2,1)=\left\{
\begin{array}{lll}
n(8nm+3m-8)/48, &  n\equiv 0\ ({\rm mod}\ 3),\ n\neq 6,9, {\rm\ and\ } \\
& m\equiv 8\ ({\rm mod}\ 16);\\
n(32nm+11m-32)/192, &  n\equiv 0\ ({\rm mod}\ 3),\ n\neq 6,9, {\rm\ and\ } \\
& m\equiv 32\ ({\rm mod}\ 64);\\
n(8nm+3m+4)/48, &  n\equiv 0\ ({\rm mod}\ 3),\ n\neq 6,9,\ m>4, \\
& m\equiv 4,20\ ({\rm mod}\ 48), {\rm\ and\ }  m/4\in S\\
\end{array}
\right.
$$
codewords, where $S$ is the set of positive integers such that for any $s\in S$, it holds that $s\equiv 1,5\ ({\rm mod}\ 12)$, and every prime divisor $p$ of $s$ satisfies $p\equiv 5\ ({\rm mod}\ 8)$, or $p\equiv 1\ ({\rm mod}\ 8)$ and $4|ord_p(2)$.
\end{Corollary}

\proof By Lemmas \ref{8 mod 16}, \ref{32 mod 64} and \ref{4 mod 16}, there is an optimal $2$-D $(3\times m,3,2,1)$-OOC with $3m/2+\Psi^e(m,3,2,1)$ codewords, where
$$\Psi^e(m,3,2,1)=\left\{
\begin{array}{lll}
(3m-8)/16, & m\equiv 8\ ({\rm mod}\ 16);\\
(11m-32)/64, & m\equiv 32\ ({\rm mod}\ 64);\\
(3m+4)/16, & m\equiv 4,20\ ({\rm mod}\ 48), m>4 {\rm\ and\ }  m/4\in S.
\end{array}
\right.
$$
Then apply Lemma \ref{n-0-mod 3} to obtain an optimal $2$-D $(n\times m,3,2,1)$-OOC with $n(nm+2\Psi^e(m,3,2,1))/6$ codewords. \qed

Combining the results of Lemmas \ref{bound}, \ref{n=2}, \ref{4} and Corollary \ref{cor:=0 mod 3}, one can complete the proof of Theorem \ref{main result}.

\section{Concluding remarks}

In Section 5.2, we present several direct constructions for optimal $2$-D $(3\times m,3,2,1)$-OOCs. Start from these 2-D OOCs and then apply Construction \ref{from 3-SCGDD} by using $m$-cyclic 3-GDDs to get some optimal $2$-D $(n\times m,3,2,1)$-OOC for any $n\equiv 0\ ({\rm mod}\ 3)$ and $n\geq 12$.

Actually to deal with the case of $n\equiv 4\ ({\rm mod}\ 6)$, by similar arguments, we can have the following lemma, in which all the input OOCs are required to attain the upper bound in Lemma \ref{0mod4}. The reader can check that the output OOC can also attain the upper bound in Lemma \ref{0mod4}.

\begin{Lemma}\label{concluding remark}
Let $m\equiv 0\ ({\rm mod}\ 4)$. Suppose that there exist
\begin{enumerate}
\item[$(1)$] an $m$-cyclic $3$-GDD of type $(6m)^{u} (4m)^{1}$, which has $2mu(1+3u)$ base blocks;
\item[$(2)$] an optimal $2$-D $(6\times m,3,2,1)$-OOC with $6m+2\Psi^e(m,3,2,1)$ codewords;
\item[$(3)$] an optimal $2$-D $(4\times m,3,2,1)$-OOC with $\lfloor (8m+4\Psi^e(m,3,2,1))/3 \rfloor$ codewords.
\end{enumerate}
Then there is an optimal $2$-D $((6u+4)\times m,3,2,1)$-OOC with $\lfloor (6u+4)((6u+4)m+2\Psi^e(m,3,2,1))/6 \rfloor$ codewords.
\end{Lemma}

By standard design theoretic techniques, it is readily checked that an $m$-cyclic $3$-GDD of type $(6m)^u(4m)^1$ exists for any $m\equiv 0\ ({\rm mod}\ 4)$ and $u\geq 3$. However, as we pointed out in Remark \ref{remark:difficult}, it seems to be difficult to find effective recursive constructions, especially filling constructions, for optimal $2$-D $(4\times m,3,2,1)$-OOCs and optimal $2$-D $(6\times m,3,2,1)$-OOCs. Even though by tedious computation and analysis, direct constructions should always work, that could not be a good way. A better technique is desired for this problem.

On the other hand, by Lemma \ref{0mod4}, we see that the size of an optimal $2$-D $(3\times m,3,2,1)$-OOC relies heavily on the size of an optimal equi-difference $1$-D $(m,3,2,1)$-OOC. The latter is closely related to optimal equi-difference CAC$(m,3)$s. This provides another motivation to study equi-difference CAC$(m,3)$s.

\end{document}